\documentclass[10pt]{article}
\usepackage{graphicx}
\title{An asymptotic expansion for the generalised quadratic Gauss sum revisited}
\author{\sc R. B.\ Paris \\
{\em School of Engineering, Computing and Applied Mathematics,} \\
{\em University of Abertay Dundee, Dundee DD1 1HG, UK}}
\begin{document}
\def\f#1#2{\mbox{${\textstyle \frac{#1}{#2}}$}}
\def\dfrac#1#2{\displaystyle{\frac{#1}{#2}}}
\def\boldal{\mbox{\boldmath $\alpha$}}
\newcommand{\bee}{\begin{equation}}
\newcommand{\ee}{\end{equation}}
\newcommand{\la}{\lambda}
\newcommand{\ka}{\kappa}
\newcommand{\al}{\alpha}
\newcommand{\ba}{\beta}
\newcommand{\om}{\omega}
\newcommand{\fr}{\frac{1}{2}}
\newcommand{\fs}{\f{1}{2}}
\newcommand{\g}{\Gamma}
\newcommand{\br}{\biggr}
\newcommand{\bl}{\biggl}
\newcommand{\ra}{\rightarrow}
\newcommand{\mbint}{\frac{1}{2\pi i}\int_{c-\infty i}^{c+\infty i}}
\newcommand{\mbcint}{\frac{1}{2\pi i}\int_C}
\newcommand{\mboint}{\frac{1}{2\pi i}\int_{-\infty i}^{\infty i}}
\newcommand{\gtwid}{\raisebox{-.8ex}{\mbox{$\stackrel{\textstyle >}{\sim}$}}}
\newcommand{\ltwid}{\raisebox{-.8ex}{\mbox{$\stackrel{\textstyle <}{\sim}$}}}
\renewcommand{\topfraction}{0.9}
\renewcommand{\bottomfraction}{0.9}
\renewcommand{\textfraction}{0.05}
\newcommand{\mcol}{\multicolumn}
\date{}
\maketitle
\pagestyle{myheadings}
\markboth{\hfill \sc R. B.\ Paris \hfill}
{\hfill \sc Asymptotics of the Gauss sum \hfill}
\begin{abstract}
An asymptotic expansion for the generalised  quadratic Gauss sum 
$$S_N(x,\theta)=\sum_{j=1}^{N} \exp (\pi ixj^2+2\pi ij\theta),$$
where $x$, $\theta$ are real and $N$ is a positive integer, is obtained as $x\ra 0$ and $N\ra\infty$ such that $Nx$ is finite. The form of this expansion holds for all values of $Nx+\theta$ and, in particular, in the neighbourhood of integer values of $Nx+\theta$. A simple bound for the remainder in the expansion is derived. Numerical results are presented to demonstrate the accuracy of the expansion and the sharpness of the bound.
 
\end{abstract}
 
\vspace{0.5cm}
 
\vspace{0.2cm} 
 

\noindent $\,$\hrulefill $\,$
 
\vspace{0.2cm}
 
\begin{center}
{\bf 1. \  Introduction}
\end{center}
\setcounter{section}{1}
\setcounter{equation}{0}
\renewcommand{\theequation}{\arabic{section}.\arabic{equation}}
We consider the asymptotic expansion of the generalised quadratic Gauss sum
\begin{equation}\label{e11}
S_N(x,\theta)=\sum_{j=1}^{N}f(j),\quad f(t):= \exp (\pi ixt^2+2\pi i\theta t),\quad 0<x<1,\ \ -\fs\leq\theta\leq\fs,
\end{equation}
where $N$ is a positive integer, as $x\ra 0$ and $N\ra\infty$, such that the quantity $Nx$ is finite.
Applications of the above exponential sum arise in
various number-theoretic contexts and in the study of disorder in dynamical systems.
 
The sum $S_N(x,\theta)$ has a long history that goes back to Gauss, who evaluated the sum 
corresponding to $x=2/N$ when $\theta=0$. 
The results of Gauss were generalised for rational $x=M/N$, where $M$ and $N$ are 
relatively prime, into the well-known Cauchy-Kronecker formula \cite{A}  
\[S_N(x,0)=\frac{e^{\pi i/4}}{\surd x} S_M\left(-\frac{1}{x},0\right)\qquad (x=M/N,\ MN\ \mbox{even}).\] 
When its terms are regarded as unit vectors in the 
complex plane, the patterns produced by the partial sums of (\ref{e11}) for fixed $x$ as $N\ra\infty$ often result in a superposition of spirals (or ``curlicues'') that can be highly intricate; see \cite{CK2, CK, DM, Le}. The scalings of this hierarchy of spirals are 
found to depend delicately on the arithmetic nature of $x$ \cite{CK}. 
When $x=p/q$, where $p$ and $q$ are relatively prime, and $\theta=0$ the trace  
is relatively simple: when $pq$ is even the spiral pattern is regular and `diffuses' in the complex 
plane away from the origin in blocks, whereas when $pq$ is odd the pattern is periodic and repeats itself indefinitely
as $N\ra\infty$. When $x$ is irrational a more complicated pattern emerges that seems to exhibit a random-walk behaviour; see \cite{CK, P}.
 
Estimates for the growth of $S_N(x,\theta)$ when $N$ is large and $x$ is fixed in the range $0<x<1$
are obtained by employing a renormalisation process based on the approximate functional relation \cite{HL}
\begin{equation}\label{e12}
S_N(x,\theta)=\frac{e^{-\pi i\theta^2/x+\pi i/4}}{\surd x}S_{\lfloor Nx\rfloor}\left(-\frac{1}{x},\frac{\theta}{x}\right)
+O\left(\frac{1+|\theta|}{\surd x}\right).
\end{equation}
This transformation shows that the sum $S_N(x,\theta)$ over $N$ terms can be approximated by a similar sum taken over 
$\lfloor Nx \rfloor$ terms with the variable $x$ replaced by $-1/x$ and $\theta$ by $\theta/x$.
Repeated application of (\ref{e12}), making use of the simple symmetry properties satisfied by
(\ref{e11}) to maintain $x$ in the interval $0<x<1$ at each stage, enables the representation of $S_N(x,\theta)$ in 
terms of a steadily decreasing number of terms. In this way it was shown in \cite{HL} that $S_N(x,\theta)=o(N)$
for any irrational $x$, with more precise order estimates depending on the detailed arithmetic structure
of $x$.

The problem that concerns us here is the asymptotic estimation of $S_N(x,\theta)$ for
$x\ra 0$ when $N\ra\infty$ such that $Nx$ is finite. An early paper dealing with estimates for $S_N(x,\theta)$
when $0<x<1$ 
is that of Fiedler {\it et al.} \cite{FJK}, and more recently that in \cite[\S 2.2]{K}, but their error terms
are too large for our purposes when $x\ra 0$. Following on from the gross estimates in
\cite{Le}, the leading terms in the expansion in the case $\theta=0$ were obtained in \cite{SZ} when $Nx<1$.
An expansion for $S_N(x,0)$ valid as $x\ra 0$  and 
finite $Nx$ was obtained in \cite{CK} and \cite[Theorem 4]{EMY}, although the remainder term was 
left as an order estimate. 

In this paper, we revisit the expansion of $S_N(x,\theta)$ as $x\ra 0$ and $N\ra\infty$ such that 
$Nx=O(1)$ obtained in \cite{P}. The sum $S_N(x,\theta)$ is expressed exactly as a series of
complementary error functions with argument proportional to $x^{-1/2}$, so that in the small-$x$ limit we may employ
the well-known asymptotics of the complementary error function in the form \cite[\S 7.12(i)]{DLMF}
\bee\label{e1a}
e^{z^2}\mbox{erfc}\,(z)= \frac{1}{\surd\pi}\sum_{r=0}^{n-1}(-)^r\frac{\g(r+\fs)}{\g(\fs)} z^{-2r-1}+{\hat T}_n(z)
\qquad (|z|\ra\infty),
\ee
where 
\bee\label{e1b}
|{\hat T}_n(z)|\leq \frac{\g(n+\fs)}{\pi}\,|z|^{-2n-1}\qquad (|\arg\,z|\leq\f{1}{4}\pi)
\ee
and $n$ is a positive integer; see also \cite[p.~111]{Olv}.
In \cite{P}, the coefficients in the resulting expansion were expressed in terms of even-order derivatives of $\cot \pi\xi$, where $\xi=Nx+\theta$, which presented a complication when
$\xi$ passes through integer values. In addition, the remainder in the expansion was not expressed as a convenient bound. Here we remedy these deficiencies and give the expansion in a form with coefficients that do not present any difficulty in computation in the neighbourhood of integer values of $\xi$.

In order to make the paper reasonably self-contained, we repeat in Section 2 the derivation of the representation
of $S_N(x,\theta)$ in terms of complementary error functions given in \cite{P}. In Section 3, we establish the central result of the paper in the following theorem.
\newtheorem{theorem}{Theorem}
\begin{theorem}
Let $S_N(x,\theta)$ be the sum defined in (\ref{e11}), where $0<x<1$, $-\fs\leq\theta\leq\fs$ and $N$ is a positive integer. Further, let $\xi:=Nx+\theta$, $M=[\xi]$ be the nearest integer part of $\xi$ and
$\epsilon:=\xi-M$, where $-\fs<\epsilon\leq\fs$. Then, as $x\rightarrow 0$ and $N\ra\infty$, such that $Nx$ is finite, we have the expansion valid for $M\geq 0$ and $n\geq 1$
\[S_N(x,\theta)-\frac{e^{-\pi i\theta^2/x+\pi i/4}}{\surd x} S_M\left(-\frac{1}{x}, \frac{\theta}{x}\right)-\frac{1}{2}(f(N)-1)\hspace{4cm}\]
\bee\label{e211}
\hspace{3.5cm}-\frac{e^{\pi i/4}}{2\surd x}\{E(\theta)-f(N)E(\epsilon)\}
=\frac{1}{2\pi i}\sum_{r=0}^{n-1}\frac{\g(r+\fs)}{\g(\fs)}\left(\frac{x}{\pi i}\right)^r C_r+R_n,
\ee
where $E(t):=e^{-\pi it^2/x}\,\mbox{erfc}\,(e^{\pi i/4} t\sqrt{\pi/x})$. The coefficients $C_r$ are given by
\bee\label{e212}
C_r=f(N)\Delta_r^-(\epsilon)-\Delta_r^-(\theta)\qquad(r\geq 0)
\ee
and the remainder $R_n$ satisfies the bound
\bee\label{e213}
|R_n|\leq \frac{(\frac{1}{2})_n}{2\pi}\left(\frac{x}{\pi}\right)^n (\Delta_n^+(\epsilon)+\Delta_n^+(\theta))\qquad(n\geq 1),
\ee
where the quantities $\Delta_r^\pm(\lambda)$ are defined in (\ref{e28}) and (\ref{e27}).
\end{theorem}
In Section 4, we present numerical results to demonstrate the accuracy of the above expansion and also the sharpness
of the bound on the remainder term $R_n$.
\vspace{0.6cm}
 
\begin{center}
{\bf 2. \ A representation for $S_N(x,\theta)$}
\end{center}
\setcounter{section}{2}
\setcounter{equation}{0}
\renewcommand{\theequation}{\arabic{section}.\arabic{equation}}
Let $\xi:=Nx+\theta$, $M=[\xi]$, $\epsilon=\xi-[\xi]$, where $[\xi]$ denotes the {\it nearest\/} integer part of $\xi$ and $-\fs<\epsilon\leq\fs$. Define also the function
\[
E(t):=e^{-\pi it^2/x}\,\mbox{erfc}\,(\om t\sqrt{\pi/x}),\qquad \om=e^{-\pi i/4},
\]
\bee\label{e20}
E(0)=1,\qquad E(-t)=2e^{-\pi it^2/x}-E(t),
\ee
where erfc is the complementary error function. The reflection formula follows from the well-known result
$\mbox{erfc}\,(z)=2-\mbox{erfc}\,(-z)$. From (\ref{e1a}), we have the expansion for $x^{-1/2}t\ra+\infty$
\bee\label{e20a}
E(t)=\frac{1}{\surd\pi}\sum_{r=0}^{n-1}(-)^r\frac{\g(r+\fs)}{\g(\fs)}\left(\frac{ix}{\pi t^2}\right)^{r+\frac{1}{2}}+T_n(t)\qquad (n=1, 2, \ldots ),
\ee
where from (\ref{e1b})
\bee\label{e20b}
|T_n(t)|\leq \frac{\g(n+\fs)}{\pi}\,\left(\frac{x}{\pi t^2}\right)^{n+\frac{1}{2}}.
\ee

An application of Cauchy's theorem shows that
\[\sum_{j=1}^{N-1}f(j)=\frac{1}{2i} \int_{\cal C}\cot (\pi t) f(t)\,dt,\]
where $f(t)$ is defined in (\ref{e11}) and ${\cal C}$ is a closed path encircling only the poles of the integrand at $t=1,2, \ldots , N-1$.
We deform
the path ${\cal C}$ into a parallelogram with two sides inclined at $\f{1}{4}\pi$ to the real axis; see \cite{P}. The vertices are situated at $\pm Pe^{\pi i/4}$, $N\pm Pe^{\pi i/4}$
($P>0$) and there are semi-circular indentations of radius $\delta<1$ around the points $t=0$ and $t=N$. Then, denoting the upper and lower halves of the contour by ${\cal C}_1$ and ${\cal C}_2$ respectively, we find following the discussion given in \cite[p.~290]{Olv} that
\[\sum_{j=1}^{N-1}f(j)=\int_\delta^{N-\delta} f(t)\,dt+\int_{{\cal C}_1}\frac{f(t)}{1-e^{-2\pi it}}dt
+\int_{{\cal C}_2}\frac{f(t)}{e^{2\pi it}-1}dt.\]
 
Now let $P\rightarrow\infty$, so that the contributions from the parts of ${\cal C}_1$ and ${\cal C}_2$
parallel to the real axis vanish on account of the exponential decay of the factor $\exp (\pi ixt^2)$, and let $\delta\rightarrow 0$. The integrals around the indentation linking $\delta e^{\pi i/4}$ with
$\delta$ and $\delta$ with $-\delta e^{\pi i/4}$ then tend to $-\f{1}{8}f(0)$ and $-\f{3}{8}f(0)$, respectively; similarly for the indentation at $t=N$ the integrals contribute $\fs f(N)$. Thus we obtain
\bee\label{e21}
S_N(x, \theta)=\sum_{j=1}^N f(j)=\fs(f(N)-1)+J_N+e^{\pi i/4}(I_N-I_0),
\ee
where the integral 
\bee\label{e21c}
J_N:=\int_0^N f(t)\,dt=\frac{e^{\pi i/4}}{2\surd x}\{E(\theta)-f(N) E(\xi)\}
\ee
and we have defined
\[I_j:=\int_0^\infty\frac{F_j(\tau)}{e^{2\pi\omega\tau}-1}d\tau\qquad (j=0,N)\]
with
\[ F_j(\tau):=f(j-\tau e^{\pi i/4})-f(j+\tau e^{\pi i/4})=2e^{-\pi x\tau^2} f(j) \sinh \{2\pi(jx+\theta)\omega\tau\}.\]

It now remains to evaluate the integrals $I_0$ and $I_N$.
If we expand the factor $(e^{2\pi\omega\tau}-1)^{-1}$ as a finite geometric series together with a remainder we find, for positive integer $K$,
\bee\label{e23a}
I_j=\sum_{k=1}^K \int_0^\infty e^{-2\pi k\omega\tau} F_j(\tau)\,d\tau + \int_0^\infty\frac{e^{-2\pi K\omega\tau}F_j(\tau)}{e^{2\pi\omega\tau}-1}d\tau.
\ee
The first term on the left-hand side of this expression becomes upon insertion of the definition of $F_j(\tau)$
\[2f(j)\sum_{k=1}^K\int_0^\infty e^{-\pi x\tau^2-2\pi k\omega\tau} \sinh \{2\pi(jx+\theta)\omega\tau\}
\,d\tau\hspace{5cm}\]
\bee\label{e23}
\hspace{5cm}=\frac{f(j)}{2\surd x} \sum_{k=1}^K\{E(k-jx-\theta)-E(k+jx+\theta)\}.
\ee

The remainder term in (\ref{e23a}) is given by
\[H_K:=\int_0^\infty\frac{e^{-2\pi K\omega\tau}F_j(\tau)}{e^{2\pi\omega\tau}-1}d\tau=f(j)\int_0^\infty
e^{-\pi\tau(x\tau+\omega)}e^{-2\pi(K-jx-\theta)\omega\tau}\,G(\tau)\,d\tau,\]
where $$G(\tau):=e^{-2\pi(jx+\theta)\omega\tau}\sinh \{2\pi(jx+\theta)\omega\tau\}/\sinh (\pi\omega\tau).$$
Now $G(0)=2(jx+\theta)$ and $G(\tau)\sim e^{-\pi\omega\tau}$ as $\tau\rightarrow+\infty$. It is also easy to see (we omit these details) that $|G(\tau)|\leq G(0)$ for $\tau\in[0,\infty)$. Then, provided $K>jx+\theta$ ($j=0, N$) it follows that
\[|H_K|< \int_0^\infty e^{-2\pi(K-jx-\theta)\omega_r\tau} |G(\tau)|\,d\tau\leq\frac{2^\frac{1}{2} G(0)}{2\pi(K-jx-\theta)}\]
where $\omega_r=1/\surd 2$, with the result that $H_K\rightarrow 0$ as $K\rightarrow\infty$. Therefore, from (\ref{e23}), we obtain
\begin{equation}\label{e24}
I_j=\frac{f(j)}{2\surd x} \sum_{k=1}^\infty \{E(k-jx-\theta)-E(k+jx+\theta)\}\qquad (j=0, N).
\end{equation}
Substitution of (\ref{e24}) with $j=0, N$  into (\ref{e21}) then yields the desired representation of $S_N(x,\theta)$ in terms of complementary error functions.

From (\ref{e20a}) we see that the terms in (\ref{e24}) are $O(k^{-2})$ as $k\ra\infty$.
\vspace{0.6cm}
 
\begin{center}
{\bf 3. \ The expansion of $S_N(x,\theta)$ as $x\ra 0$ with $Nx$ finite}
\end{center}
\setcounter{section}{3}
\setcounter{equation}{0}
\renewcommand{\theequation}{\arabic{section}.\arabic{equation}}
We define the quantities for $|\lambda|<1$
\bee\label{e28}
\Delta_r^+(\lambda):=\zeta(2r+1, 1+\lambda)+\zeta(2r+1, 1-\lambda)\quad (r\geq 1),
\ee
and
\bee\label{e27}
\Delta_r^-(\lambda):=\left\{\begin{array}{ll} \pi \cot \pi\lambda-\lambda^{-1} &  (r=0)\\
\\
\zeta(2r+1, 1+\lambda)-\zeta(2r+1, 1-\lambda) & (r\geq 1),\end{array}\right.
\ee
where $\zeta(s,a)=\sum_{k=0}^\infty(k+a)^{-s}$ ($\Re (s)>1$) is the Hurwitz zeta function. Note that
$\Delta_r^-(0)=0$ for $r\geq 0$ and $\Delta_r^+(0)=2\zeta(2r+1)$, where $\zeta(s)$ is the Riemann zeta function.

When $j=0$, we have from (\ref{e24})
\[I_0=\frac{1}{2\surd x}\sum_{k=1}^\infty \{E(k-\theta)-E(k+\theta)\}.\]
In the limit $x\rightarrow 0$, the arguments of the complementary error functions contained in $E(k\pm\theta)$ have large modulus for $k\geq 1$ and phase equal to $-\f{1}{4}\pi$, since $-\fs\leq\theta\leq\fs$. 
Employing the expansion (\ref{e20a}), we then obtain
\begin{equation}\label{e25}
I_0=-\frac{e^{\pi i/4}}{2\pi}\sum_{r=1}^{n-1} \frac{\g(r+\fs)}{\g(\fs)} \left(\frac{x}{\pi i}\right)^r c_r(\theta)+{\cal R}_n(\theta)\qquad (x\ra 0),
\end{equation}
where
\begin{eqnarray}
c_r(\theta)=\mathop{\sum_{k=-\infty}^\infty}_{\scriptstyle k\neq 0} (k+\theta)^{-2r-1}
&=&\sum_{k=0}^\infty (k+1+\theta)^{-2r-1}-\sum_{k=0}^\infty (k+1-\theta)^{-2r-1}\nonumber\\
&=&\Delta_r^-(\theta).\label{e25a}
\end{eqnarray}
In the case $r=0$ the sums must be interpreted in the principal value sense $\lim_{s\rightarrow\infty}\sum_{k=-s}^s a_k$ to yield the evaluation $c_0(\theta)=\pi \cot \pi\theta-1/\theta$.
The remainder term ${\cal R}_n(\theta)$ is given by
\[{\cal R}_n(\theta)=\frac{1}{2\surd x}\sum_{k=1}^\infty\{T_n(k-\theta)-
T_n(k+\theta)\}\]
and, from (\ref{e20b}), therefore satisfies the bound
\bee\label{e25b}
|{\cal R}_n(\theta)|
\leq\frac{(\frac{1}{2})_n}{2\pi}\,\left(\frac{x}{\pi}\right)^n\mathop{\sum_{k=-\infty}^\infty}_{\scriptstyle k\neq 0} |k+\theta|^{-2n-1}=\frac{(\frac{1}{2})_n}{2\pi}\,\left(\frac{x}{\pi}\right)^n\,\Delta_n^+(\theta),
\ee
where $(a)_r=\g(a+r)/\g(a)$ is the Pochhammer symbol.
 
Proceeding in a similar manner when $j=N$, we have
\begin{eqnarray*}
I_N&=&\frac{f(N)}{2\surd x}\sum_{k=1}^\infty \{E(k-\xi)-E(k+\xi)\}\\
&=&\frac{f(N)}{2\surd x}\left\{\sum_{k=1}^M\{2e^{-\pi i(k-\xi)^2/x}-E(\xi-k)\}+\sum_{k=M+1}^\infty E(k-\xi)-\sum_{k=1}^\infty E(k+\xi)\right\}.
\end{eqnarray*}
Here we have made use of the reflection formula in (\ref{e20}) to separate off the error functions in $E(k-\xi)$ corresponding to $k\leq M$ (when $M\geq 1$). Upon noting that
\[f(N) e^{-\pi i(k-\xi)^2/x}=e^{-\pi i\theta^2/x}\,e^{-\pi ik^2/x+2\pi ki\theta/x},\]
we obtain when $M\geq 1$
\bee\label{e26b}
I_N=\frac{e^{-\pi i\theta^2/x}}{\surd x} S_M\left(-\frac{1}{x}, \frac{\theta}{x}\right)
-\frac{f(N)}{2\surd x}\left\{\sum_{k=1}^\infty E(k+\xi)+\sum_{k=1}^{M}E(\xi-k)-\sum_{k=M+1}^\infty E(k-\xi)\right\}.
\end{equation}
If we now extract from the second sum in curly braces in (\ref{e26b}) the error function $E(\xi-k)$ corresponding to $k=M$ (that is, $E(\epsilon)$) and use the evaluation of the integral $J_N$ in (\ref{e21c}), we can write
\[e^{\pi i/4}I_N+J_N=\frac{e^{-\pi i\theta^2/x+\pi i/4}}{\surd x} S_M\left(-\frac{1}{x}, \frac{\theta}{x}\right)+\frac{e^{\pi i/4}}{2\surd x}\{E(\theta)-f(N)E(\epsilon)\}\hspace{2cm}\]
\begin{equation}\label{e26}
\hspace{2cm}-\frac{e^{\pi i/4}f(N)}{2\surd x}\left\{\sum_{k=0}^\infty E(k+\xi)+\sum_{k=1}^{M-1}E(\xi-k)-\sum_{k=M+1}^\infty E(k-\xi)\right\},
\end{equation}
where the term involving $f(N) E(\xi)$ from $J_N$ has been absorbed into the first sum in curly braces.
 
Then in a similar manner to the determination of the expansion of $I_0$ in (\ref{e25}) we find
\[\frac{1}{2\surd x}\left\{\sum_{k=0}^\infty E(k+\xi)+\sum_{k=1}^{M-1}E(\xi-k)-\sum_{k=M+1}^\infty E(k-\xi)\right\}\hspace{4cm}\]
\bee\label{e25c}
\hspace{2cm}=\frac{e^{\pi i/4}}{2\pi}\sum_{r=0}^{n-1} \frac{\g(r+\fs)}{\g(\fs)}\left(\frac{x}{\pi i}\right)^r c_r(\epsilon)+{\cal R}_n(\epsilon),
\ee
where, recalling that $\xi=M+\epsilon$, $M=[\xi]$,
\bee\label{e25d}
c_r(\epsilon)=\mathop{\sum_{k=-\infty}^\infty}_{\scriptstyle k\neq -M}(k+\xi)^{-2r-1}
=\mathop{\sum_{k=-\infty}^\infty}_{\scriptstyle k\neq 0}(k+\epsilon)^{-2r-1}=\Delta_r^-(\epsilon)
\ee
and the remainder ${\cal R}_n(\epsilon)$ satisfies the bound
\begin{eqnarray}
|{\cal R}_n(\epsilon)|&\leq&\frac{(\frac{1}{2})_n}{2\pi}\left(\frac{x}{\pi}\right)^n
\mathop{\sum_{k=-\infty}^\infty}_{\scriptstyle k\neq -M}|k+\xi|^{-2n-1}
=\frac{(\frac{1}{2})_n}{2\pi}\left(\frac{x}{\pi}\right)^n
\mathop{\sum_{k=-\infty}^\infty}_{\scriptstyle k\neq 0}|k+\epsilon|^{-2n-1}\nonumber\\
&=&\frac{(\frac{1}{2})_n}{2\pi}\left(\frac{x}{\pi}\right)^n\,\Delta_n^+(\epsilon).\label{e25e}
\end{eqnarray}

Combination of (\ref{e26}) and (\ref{e25c}) then yields the expansion when $M\geq 1$
\[e^{\pi i/4}I_N+J_N=\frac{e^{-\pi i\theta^2/x+\pi i/4}}{\surd x} S_M\left(-\frac{1}{x}, \frac{\theta}{x}\right)+\frac{e^{\pi i/4}}{2\surd x}\{E(\theta)-f(N)E(\epsilon)\}\]
\bee\label{e215}
+\frac{1}{2\pi i}\sum_{r=0}^{n-1} \frac{\g(r+\fs)}{\g(\fs)} \left(\frac{x}{\pi i}\right)^r c_r(\epsilon)
+{\cal R}_n(\epsilon).
\ee
In the case $M=0$ (when $\xi=\epsilon$), the sum $S_M\equiv 0$ and from (\ref{e26b}) we have
\[e^{\pi i/4}I_N+J_N=\frac{e^{\pi i/4}}{2\surd x}\{E(\theta)-f(N)E(\epsilon)\}
-\frac{e^{\pi i/4}f(N)}{2\surd x}\left\{\sum_{k=1}^\infty E(k+\epsilon)-\sum_{k=1}^\infty E(k-\epsilon)\right\}.\]
It is easily seen that we obtain the same expansion as (\ref{e215}).

The form of the coefficients in (\ref{e25a}) and (\ref{e25d}) with $r\geq 1$ presents no difficulty in computation in the neighbourhood of integer values of $\xi$ where $\epsilon \simeq 0$, in contrast to those given in \cite{P} which involved even derivatives of $\cot \pi\xi$. Although the coefficients $c_0(\epsilon)$ and $c_0(\theta)$ have a removable singularity at $\epsilon=0$ and $\theta=0$ their computation is straightforward. 

If we now define the coefficients $C_r$ and the remainder $R_n$ by
\[C_r:=f(N)c_r(\epsilon)-c_r(\theta),\qquad R_n:=e^{\pi i/4}\{f(N) {\cal R}_n(\epsilon)-{\cal R}_n(\theta)\},\] 
then we see that
\bee\label{e214}
C_r=f(N)\Delta_r^-(\epsilon)-\Delta_r^-(\theta) \quad (r\geq 0)
\ee
and
\bee\label{214a}
|R_n|\leq |{\cal R}_n(\epsilon)|+|{\cal R}_n(\theta)|\leq \frac{(\fs)_n}{2\pi} \left(\frac{x}{\pi}\right)^n\{\Delta_n^+(\epsilon)+\Delta_n^+(\theta)\}\quad (n\geq 1).
\ee
Combination of (\ref{e21}), (\ref{e25}) and (\ref{e215}), together with the above definitions of $C_r$ and the bound
on $R_n$, then gives the expansion of $S_N(x,\theta)$ stated in 
Theorem 1. We remark that the terms $E(\theta)$ and $E(\epsilon)$ have been left unexpanded as $x\ra 0$ in (\ref{e215})
and in Theorem 1, since for small values of $\theta$ and $\epsilon=o(x^{1/2})$ these quantities can no longer be approximated by (\ref{e20a}).
 
\vspace{0.3cm}

\begin{center}
{\bf 4. \ Numerical results and discussion}
\end{center}
\setcounter{section}{4}
\setcounter{equation}{0}
\renewcommand{\theequation}{\arabic{section}.\arabic{equation}}
In order to demonstrate the accuracy of the expansion in Theorem 1, we define the quantity ${\cal S}$ by
\bee\label{e41}
{\cal S}:=S_N(x,\theta)-\frac{e^{-\pi i\theta^2/x+\pi i/4}}{\surd x} S_M\left(-\frac{1}{x}, \frac{\theta}{x}\right)-\frac{1}{2}(f(N)-1)
-\frac{e^{\pi i/4}}{2\surd x}\{E(\theta)-f(N) E(\epsilon)\}.
\ee
Then from Theorem 1 we have the expansion as $x\ra 0$ and $N\ra\infty$ such that $Nx$ is finite
\bee\label{e42}
{\cal S}=\frac{1}{2\pi i}\sum_{r=0}^{n-1}\frac{\g(r+\fs)}{\g(\fs)}\left(\frac{x}{\pi i}\right)^r C_r+R_n,
\ee
where the coefficients $C_r$ are defined in (\ref{e212}) and the remainder $R_n$ satisfies the bound in (\ref{e213}).
We remark that the bound in (\ref{e213}) is explicitly independent of $N$.
In Table 1, we show the absolute value of the error in the computation of ${\cal S}$ using the expansion (\ref{e42}) truncated after $n$ terms for two different sets of values of $x$, $\theta$, summation index $N$ and different levels $n$. The exact value of $S_N(x,\theta)$ was obtained by high-precision summation of (\ref{e11}).  
In Table 2, we compare the absolute values of the remainder $R_n$ calculated from (\ref{e42})
and its bound to illustrate the sharpness of (\ref{e213}).
\begin{table}[th]
\caption{\footnotesize{Values of the absolute error in the computation of ${\cal S}$ by 
(\ref{e42}) for different truncation index $n$.}}
\begin{center}
\begin{tabular}{r|c|c|c}
\hline
\mcol{1}{c|}{} & \mcol{1}{c|}{$x=1/(250\surd\pi)$} & \mcol{1}{c|}{$x=1/(250\surd\pi)$} & \mcol{1}{c}{$x=1/(250\surd 3)$}
\\
\mcol{1}{c|}{} & \mcol{1}{c|}{$N=7300$,\ \ $\theta=-0.125$} & \mcol{1}{c|}{$N=7430$,\ \ $\theta=0.25$}
& \mcol{1}{c}{$N=6000,\ \ \theta=0$}
\\
\mcol{1}{c|}{$n$} & \mcol{1}{c|}{$\xi\doteq 16.349$} & \mcol{1}{c|}{$\xi\doteq 17.018$} & \mcol{1}{c}{$\xi\doteq 6.928$}
\\
[.15cm]\hline
&&&\\[-0.25cm]
1  & $2.216\times 10^{-4}$  & $1.198\times 10^{-4}$ & $1.386\times 10^{-5}$ \\
2  & $5.642\times 10^{-7}$  & $2.527\times 10^{-7}$ & $1.221\times 10^{-8}$  \\
3  & $2.346\times 10^{-9}$  & $\ \,8.332\times 10^{-10}$ & $\ \,1.590\times 10^{-11}$ \\
4  & $\ \,1.369\times 10^{-11}$ & $\ \,3.752\times 10^{-12}$ & $\ \,2.708\times 10^{-14}$ \\
6  & $\ \,9.569\times 10^{-16}$ & $\ \,1.509\times 10^{-16}$ & $\ \,1.420\times 10^{-19}$ \\
8  & $\ \,1.334\times 10^{-19}$ & $\ \,1.194\times 10^{-20}$ & $\ \,1.360\times 10^{-24}$ \\
10 & $\ \,3.096\times 10^{-23}$ & $\ \,1.568\times 10^{-24}$ & $\ \,2.082\times 10^{-29}$ \\
[.15cm]\hline
\end{tabular}
\end{center}
\end{table}
\begin{table}[th]
\caption{\footnotesize{The absolute values of $R_n$ and the bound in 
(\ref{e213}) for different truncation index $n$.}}
\begin{center}
\begin{tabular}{r|ll|ll}
\hline
\mcol{1}{c|}{} & \mcol{2}{c|}{$x=1/(250\surd\pi)$,\ \ $\theta=-0.125$} & \mcol{2}{c}{$x=1/(250\surd\pi)$,\ \ $\theta=0.25$}
\\
\mcol{1}{c|}{$$} & \mcol{2}{c|}{$N=7300$,\ \ $\xi\doteq 16.349$} & \mcol{2}{c}{$N=7430$,\ \ $\xi\doteq 17.018$}
\\
\mcol{1}{c|}{$n$} & \mcol{1}{c}{$|R_n|$} & \mcol{1}{c|}{Bound} & \mcol{1}{c}{$|R_n|$} & \mcol{1}{c}{Bound}\\
[.05cm]\hline
&&\\[-0.25cm]
1  & $2.216\times 10^{-4}$   & $4.062\times 10^{-4}$ &  $1.200\times 10^{-4}$ &  $3.272\times 10^{-4}$ \\
2  & $5.642\times 10^{-7}$   & $7.077\times 10^{-7}$ &  $2.527\times 10^{-7}$ &  $4.137\times 10^{-7}$ \\
4  & $1.369\times 10^{-11}$  & $1.435\times 10^{-11}$ & $3.752\times 10^{-12}$ & $4.309\times 10^{-12}$ \\
6  & $9.569\times 10^{-16}$  & $9.691\times 10^{-16}$ & $1.509\times 10^{-16}$ & $1.570\times 10^{-16}$ \\
8  & $1.334\times 10^{-19}$  & $1.339\times 10^{-19}$ & $1.194\times 10^{-20}$ & $1.208\times 10^{-20}$ \\
10 & $3.096\times 10^{-23}$  & $3.100\times 10^{-23}$ & $1.568\times 10^{-24}$ & $1.574\times 10^{-24}$ \\
[.15cm]\hline
\end{tabular}
\end{center}
\end{table}
 
In the case of the classical quadratic Gauss sum ($\theta=0$), we have $\xi=Nx=M+\epsilon$ with $-\fs<\epsilon\leq\fs$ (when $M\geq 1$). From (\ref{e211}) as $x\ra 0$, $N\ra\infty$ such that $Nx$ is finite, we obtain the expansion
\[S_N(x,0)=\frac{e^{\pi i/4}}{\surd x} S_M\left(-\frac{1}{x}, 0\right)+\frac{1}{2}(f(N)-1)+\frac{e^{\pi i/4}}{2\surd x}\{1-f(N)E(\epsilon)\}\hspace{2cm}\]
\bee\label{e43}
\hspace{6cm}
+\frac{f(N)}{2\pi i}\sum_{r=0}^{n-1}\frac{\g(r+\fs)}{\g(\fs)}\left(\frac{x}{\pi i}\right)^r c_r(\epsilon)+R_n',
\ee
where $f(N)=\exp (\pi ixN^2)$. From (\ref{e213}) and the fact that $I_0\equiv 0$ when $\theta=0$, we have
\bee\label{e44}
|R_n'|\leq\frac{(\fs)_n}{2\pi}\left(\frac{x}{\pi}\right)^n\,\Delta_n^+(\epsilon).
\ee
We emphasise that the expansion in (\ref{e43}) holds for all finite values of $Nx$; see Table 1.
When $M=[\xi]=0$ and $0<\epsilon\leq\fs$ --- that is, when $Nx<\fs$ --- the sum $S_M=0$. If, in addition, $E(\epsilon)$ in (\ref{e43}) is expanded by means of (\ref{e20a}) then we obtain an expansion equivalent to that in \cite[Theorem 4]{EMY}, albeit with a bound for the remainder rather than an order estimate and coefficients $c_r$ expressed in a different form. However, the expansion of $E(\epsilon)$ by (\ref{e20a}) is only applicable when $\epsilon\gg x^{1/2}$; that is, when $N\gg x^{-1/2}$.
 
Finally, we remark that since $C_r\sim (1-|\epsilon|)^{-2r-1}$ for $r\gg 1$ (when $\epsilon$ is bounded away from zero), the optimal truncation index $r_0$ of the sum in (\ref{e42}) (corresponding to truncation at, or near, the term of least magnitude) is given by $r_0\simeq \pi(1-|\epsilon|)^2/x$. This shows that the values of the truncation index $n$ in Table 1 are highly sub-optimal and also gives an indication of the enormous accuracy that could be obtained from the expansion (\ref{e42}).

\end{document}